[1994/06/01]
\documentclass[12pt]{amsart}

        [1994/05/12 v2.2b Standard AMS font definitions]
\DeclareFontFamily{U}{msa}{}
\DeclareFontShape{U}{msa}{m}{n}
  { <5> <6> <7> <8> <9> gen * msam
    <10> <10.95> <12> <14.4> <17.28> <20.74> <24.88> msam10}{}
        [1994/05/12 v2.2b Standard AMS font definitions]
\DeclareFontFamily{U}{msb}{}
\DeclareFontShape{U}{msb}{m}{n}
  { <5> <6> <7> <8> <9> gen * msbm
    <10> <10.95> <12> <14.4> <17.28> <20.74> <24.88> msbm10}{}

\newtheorem{prop}[subsection]{Proposition}
\newtheorem{rem}[subsection]{Remark}
\newtheorem{thm}[subsection]{Theorem}

\newtheorem{firstthm}{Theorem 1}

\theoremstyle{definition}

\numberwithin{equation}{section}

\title[Automorphisms of Torelli Groups]{Automorphisms of Torelli Groups}

\author[J. D. McCarthy]{John D. McCarthy}
\address{Department of Mathematics, Michigan State University, East Lansing, MI
48824-1027}
\email{mccarthy@@math.msu.edu}

\author[W. R. Vautaw]{William R. Vautaw}
\address{Department of Mathematics \ Southeastern Louisiana University \
Hammond, LA 70402}
\email{wvautaw@selu.edu}

\date{November 13, 2003}

\keywords{mapping class groups, Torelli groups}
\subjclass{Primary 32G15; Secondary 20F38, 30F10, 57M99}

\newcommand{\abstracttext}{In this paper, we prove that each automorphism of the Torelli group
of a surface is induced by a diffeomorphism of the surface, provided that the surface is a
closed, connected, orientable surface of genus at least 3. This result was previously
announced by Benson Farb for genus at least 4 and has recently been announced by Benson Farb
and Nikolai Ivanov for  subgroups of finite index in the Torelli group for surfaces of genus
at least 5. This result is also directly analogous to previous results established for
classical braid groups, mapping class groups, and surface braid groups.}

\begin{document}

\begin{abstract}  \abstracttext  \end{abstract}

\maketitle

\section{Introduction}
\label{sec:introduction}

Let $S$ be a closed, connected, oriented surface of genus $g$. The extended mapping class group
$\mathcal{M}^*$ of $S$ is defined to be the group of isotopy classes of homeomorphisms $S
\rightarrow S$.  The mapping class group $\mathcal{M}$ of $S$ is defined to be the subgroup
of index $2$ in $\mathcal{M}$ of isotopy classes of orientation preserving homeomorphisms
$S \rightarrow S$. The Torelli group $\mathcal{T}$ of $S$ is the subgroup of $\mathcal{M}$
consisting of the isotopy classes of those orientation preserving homeomorphisms $S \rightarrow
S$ that induce the identity permutation of the first homology group of $S$. 

In this paper, we show that any given automorphism $\Psi: \mathcal{T} \rightarrow \mathcal{T}$
is induced by a homeomorphism $S \rightarrow S$ provided that $g > 2$. This result is
trivially true for $g < 2$, as $\mathcal{T}$ is trivial for $g < 2$. On the other hand, by a result
of Mess, this result is false for $g = 2$. Indeed, Mess showed that $\mathcal{T}$ is a
nonabelian free group of infinite rank for $g = 2$ \cite{me}. It follows that the
automorphism group of $\mathcal{T}$ for $g = 2$ is uncountable as it contains a copy of the
permutation group of an infinite set, any free basis for $\mathcal{T}$. An  automorphism of
$\mathcal{T}$ induced by a homeomorphism $h: S \rightarrow S$ only depends upon the isotopy
class of $h$ in $\mathcal{M}^*$. Since $\mathcal{M}^*$ is finitely generated \cite{b} and,
hence, countable, it follows that only countably many elements of the automorphism group of
$\mathcal{T}$ are induced by homeomorphisms $S \rightarrow S$.  

In several formal and informal announcements made between October 2001 and March 2002, Benson Farb
announced that he had proved that any given automorphism $\Psi: \mathcal{T} \rightarrow
\mathcal{T}$ is induced by a homeomorphism $S \rightarrow S$ provided that $g > 3$
\cite{f}. In his thesis, \cite{v1}, the second author laid the foundation for proving that it is
true for $g > 2$. This involves three basic steps, two of which were finished in \cite{v1}. In the
first step, completed in \cite{v1}, Vautaw characterized algebraically certain elements of the
Torelli group, namely powers of Dehn twists about separating curves and powers of bounding pair
maps. The second step, left unfinished in \cite{v1}, is to show that $\Psi$ induces an
automorphism $\Psi_* : \mathcal{C} \rightarrow \mathcal{C}$, where $\mathcal{C}$ is the complex of
curves of $S$. In the last step, completed in \cite{v1}, Vautaw used a theorem of Ivanov which
states that $\Psi_*$ is induced by a homeomorphism $H : S \rightarrow S$ of the surface $S$
\cite{iv1}, and concluded, under the assumption that it is possible to complete the second step,
that the automorphism of the Torelli group induced by $H$ agrees with the given automorphism
$\Psi$ (i.e. $\Psi([F]) = [H \circ F \circ H^{-1}]$ for every mapping class $[F]$ in
$\mathcal{T}$). 

The purpose of this paper is to complete the second step and, hence, establish the main
theorem, Conjecture $1$ of \cite{v1}: 

\begin{firstthm} Let $S$ be a closed, connected, orientable surface of genus $g \neq 2$, and let
$\Psi : \mathcal{T} \rightarrow \mathcal{T}$ be an automorphism of the Torelli group $\mathcal{T}$
of $S$. Then $\Psi$ is induced by an homeomorphism $h : S \rightarrow S$. That is, there exists an
homeomorphism $H : S \rightarrow S$ such that for any mapping class $[F]$ in $\mathcal{T}$, we
have $\Psi([F]) = [H \circ F \circ H^{-1}]$. 
\label{firstthm:mainthm} \end{firstthm} 

This result is analogous to results previously established for classical braid groups \cite{dg},
mapping class groups \cite{i}, \cite{iv2}, \cite{ivm}, \cite{k}, \cite{m} and surface braid groups
\cite{iivm}. As noted above, it fails for $g = 2$. Recently, Farb and Ivanov have established the
analogue of this result for subgroups of finite index in the Torelli group of surfaces of genus
at least $5$ \cite{fiv}. They obtain this result as a consequence of an analogue of Ivanov's result
on automorphisms of complex of curves \cite{iv1} (see also \cite{i} and \cite{k}). This result of
Farb and Ivanov subsumes Farb's previous result except for surfaces of genus $4$. In this work,
they use a characterization of powers of Dehn twists about separating curves and powers of bounding
pair maps (Proposition 8 of \cite{fiv}) equivalent to the second author's characterization 
(Theorem 3.5 of \cite{v1}).

\section{Preliminaries} 
\label{sec:prelims} 

We have the following results from the second author's thesis \cite{v1}. 

\begin{prop} Let $\mathcal{V}_{sep}$ denote the set of isotopy classes of
essential separating circles on $S$. Let $a$ be an element of
$\mathcal{V}_{sep}$. Then there exists a unique element $c$ of
$\mathcal{V}_{sep}$ such that either $\Psi(D_a) = D_c$ or $\Psi(D_a) = D_c^{-1}$.
\label{prop:presseptwists} \end{prop}

\begin{prop} Let $\mathcal{V}_{nonsep}$ denote the set of isotopy classes of 
nonseparating circles on $S$. Let $\mathcal{BP}$ be the set of ordered pairs
$(a,b)$ of elements of $\mathcal{V}_{sep}$ such that $\{a, b\}$ is a bounding
pair on $S$. Let $(a,b)$ be an element of $\mathcal{BP}$. Then there exists a
unique element
$(e,f)$ of $\mathcal{BP}$ such that $\Psi(D_a \circ D_b^{-1}) = D_e \circ
D_f^{-1}$.
\label{prop:presbpmaps} \end{prop}

As explained in \cite{v1}, these results are derived from the characterizations of elementary
twists in the second author's thesis, Theorems 3.5, 3.7, and 3.8 of \cite{v1}. These
characterizations are analogous to corresponding characterizations used in previous works on
automorphisms of subgroups of the mapping class group (e.g. \cite{i}, \cite{iv2}, \cite{ivm},
\cite{k}, \cite{m}). As is true for these previous characterizations, these characterizations are
established by using the theory of abelian subgroups of mapping class groups developed in the
work of Birman-Lubotzky-McCarthy \cite{blm} from Thurston's theory of surface diffeomorphisms \cite{flp}. 

\section{Two-holed tori} 
\label{sec:twoholedtori}

We shall use the following generalization of the Centralizer Lemma of \cite{v2}.  

\begin{prop} Let $S$ be a connected, closed, orientable surface of genus $g > 2$.
Let $R$ be an embedded two-holed torus in $S$ such that both components,
$a$ and $b$, of the boundary $\partial R$ of $R$ are essential separating curves
on $S$. Let $\Psi : \mathcal{T} \rightarrow \mathcal{T}$ be an automorphism of
the Torelli group $\mathcal{T}$ of $S$. Let $\Gamma_R$ denote the subgroup of
$\mathcal{T}$ consisting of mapping classes of homeomorphisms $S \rightarrow S$
which are supported on $R$. Then the restriction of $\Psi$ to $\Gamma_R$ is
induced by a homeomorphism $H : S \rightarrow S$. That is, there exists an 
homeomorphism $H : S \rightarrow S$ such that for any mapping class $[F]$ in
$\Gamma_R$, we have $\Psi([F]) = [H \circ F \circ H^{-1}]$.
\label{prop:twoholedtori} \end{prop}

The proof of this proposition follows the lines of the proof of the Centralizer Lemma of \cite{v2}. From
Propositions \ref{prop:presseptwists} and \ref{prop:presbpmaps}, it follows that we may assume without
loss of generality that $\Psi : \mathcal{T} \rightarrow \mathcal{T}$ restricts to an automorphism 
$\Psi : \Gamma_R \rightarrow \Gamma_R$. The result is established by relating the restriction $\Psi :
\Gamma_R \rightarrow \Gamma_R$ to an automorphism $\chi : \pi_1(T,x) \rightarrow \pi_1(T,x)$ of the
fundamental group $\pi_1(T,x)$ of a one-holed torus $T$. The one-holed torus $T$ corresponds to the
quotient of $R$ obtained by collapsing one of the boundary components $a$ of $R$ to a basepoint $x$ for
$T$. The relevant relationship between $\Gamma_R$ and $\pi_1(T,x)$ is obtained from the long-exact
homotopy sequence of an appropriate fibration. One of the key ingredients in the argument establishing
this relationship is the well-known contractibility of the identity component of the group of
diffeomorphisms of surfaces of negative euler characteristic \cite{ee}. It is a well-known classical fact
that an automorphism of the fundamental group $\pi_1(T,x)$ of a one-holed torus $T$ is induced by a
diffeomorphism $(T,x)$ if and only if it preserves the peripheral structure of $\pi_1(T,x)$ \cite{zvoc}. As
in the proof of the Centralizer Lemma of \cite{v2}, we show that this peripheral structure corresponds to
Dehn twists about separating curves of genus 1 in the two-holed torus $R$. Together with Propositions 
\ref{prop:presseptwists} and \ref{prop:presbpmaps}, this allows us to appeal to the aforementioned
classical fact to obtain a homeomorphism $H : (S,R,a,b) \rightarrow (S,R,a,b)$ whose restriction 
$H : (R, a, b) \rightarrow (R, a, b)$ descends to a homeomorphism $(T,x) \rightarrow (T,x)$ inducing the
automorphism $\chi : \pi_1(T,x) \rightarrow \pi_1(T,x)$. We then show, as in the proof of the Centralizer
Lemma of \cite{v2}, using Propositions \ref{prop:presseptwists} and \ref{prop:presbpmaps}, that $H$ 
satisfies the conclusion on Proposition \ref{prop:twoholedtori}. 

We note that a similar argument relating certain subgroups of mapping class groups to fundamental groups
of punctured surfaces appears in the work of Irmak-Ivanov-McCarthy on automorphisms of surface braid
groups, in which they prove the analogue of Theorem \ref{firstthm:mainthm} for certain surface braid
groups \cite{iivm}.

\section{Orientation type} 
\label{sec:ortype} 

In this section, we show how to ascribe an orientation type, $\epsilon$, to any
given automorphism $\Psi : \mathcal{T} \rightarrow \mathcal{T}$. This is
achieved by studying the action of $\Psi$ on the separating twists in
$\mathcal{T}$. We show that either (i) $\Psi$ sends right separating twists to
right separating twists or (ii) $\Psi$ sends right separating twists to left
separating twists. If (i) holds, we say that $\Psi$ is orientation-preserving
and let $\epsilon = 1$. If (ii) holds, we say that $\Psi$ is
orientation-reversing and set $\epsilon = -1$.   

The proof of this result uses the following connectivity result.

\begin{prop} Let $\mathcal{C}_{sep}$ be the subcomplex of $\mathcal{C}$ whose
simplices are those simplices of $\mathcal{C}$ all of whose vertices are isotopy
classes of essential separating circles on $S$. Then $\mathcal{C}_{sep}$ is
connected. 
\label{prop:sepconnected} \end{prop} 

\begin{proof} Let $a$ and $b$ be elements of $\mathcal{C}_{sep}$. Choose circles 
$A$ and $B$ in the isotopy classes $a$ and $b$ so that $A$ and $B$ intersect 
transversely and the geometric intersection $i(a,b)$ of $a$ and $b$ is equal 
to the number of points of intersection of $A$ and $B$. Since $A$ is a separating 
circle on $S$ and the genus of $S$ is at least $3$, $A$ separates $S$ into two 
one-holed surfaces with boundary $A$, at least one of which has genus at least $2$.
Let $R$ be one of these two-holed surfaces with boundary $A$ and genus at least $2$. 

Our goal is to find a path in $\mathcal{C}_sep$ from $a$ to $b$. Since $A$ and $B$ are separating circles, $i(a,b)$
is even. We shall prove the result by induction on $i(a,b)$. We may assume that $i(a,b) > 0$. Simple arguments show
that an appropriate path exists if $i(a,b)$ is equal to $2$ or $4$. Hence, it suffices to find a separating circle
$C$ such that $i(a,c) \leq 4$ and $i(c,b) < i(a,b)$. 

Roughly speaking, such a circle $C$ can be constructed as follows. Surger $A$ into two essential simple closed
curves $A_1$ and $A_2$ in the interior of $R$ by using one of the component arcs $J$ of $B \cap R$, so that $A \cup
A_1 \cup A_2$ is the boundary of a tubular neighborhood $N$ of $A \cup J$ in $R$. If either $A_1$ or $A_2$ is
separating, then we choose $C$ as either $A_1$ or $A_2$, whichever is separating. Otherwise, by tubing $A_1$ to
$A_2$ along an appropriately chosen arc $K$ joining $A_1$ to $A_2$ and intersecting $A$ twice, we may construct an
essential separating circle $C$ such that $i(a,c) \leq 4$ and $i(c,b) < i(a,b)$. 
\end{proof} 

\begin{rem} Farb and Ivanov announced this result in their recent paper (see Theorem 
4 of \cite{fiv}). They point out that they deduce this result as a consequence of a 
stronger result regarding separating circles of genus $1$. We note that this stronger result 
follows from Proposition \ref{prop:sepconnected}. Indeed, suppose that $a_1$,...,$a_{n + 1}$ is 
a sequence of isotopy classes of essential separating circles, $a_1$ has genus $1$, $i(a_j, a_{j + 1}) = 0$ for $1
\leq j \leq n$ and $a_{n + 1}$ has genus $1$. Suppose that $a_k$ has genus $> 1$ for some integer $k$ with $1 < k
< n + 1$. Represent $a_{k - 1}$, $a_k$ and $a_{k + 1}$ by essential separating circles $A_{k - 1}$, $A_k$ and
$A_{k + 1}$ on $S$ with $A_{k - 1}$ and $A_{k + 1}$ disjoint from $A_k$. Note that $A_k$ separates $S$ into 
two one-holed surfaces, $P$ and $Q$, each of genus $> 1$. We may assume, without loss of generality, 
that $A_{k - 1}$ is in the interior of $P$. Either (i) $A_{k + 1}$ is in the interior of $P$ or (ii) $A_{k + 1}$ is
in the interior of $Q$. In case (i), we may replace $a_k$ by the isotopy class of some separating circle of
genus $1$ in $Q$.  In case (ii), we may delete $a_k$. In either case, we obtain a path in $\mathcal{C}_{sep}$ from
$a_1$ to $a_{n + 1}$ with fewer vertices of genus $> 1$ than the chosen path $a_1$,...,$a_{n + 1}$. Hence, the
stronger result follows by induction on the number of vertices of genus $> 1$ in the chosen path.
\end{rem}  
 
Using Proposition \ref{prop:sepconnected}, we may now establish the main result
of this section.

\begin{prop} Let $\Psi : \mathcal{T} \rightarrow \mathcal{T}$ be an automorphism
of the Torelli group $\mathcal{T}$. Then there exists an
$\epsilon$ in $\{ -1, 1\}$ such that for each essential separating circle
$a$ on $S$, $\Psi(D_a) = D_c^\epsilon$ for some essential separating circle $c$
on $S$. 
\label{prop:ortype} \end{prop}  

\begin{proof} Let $a$ be an essential separating circle on $S$. By Proposition
\ref{prop:presseptwists}, either $\Psi(D_a) = D_c$ or $\Psi(D_a) = D_c^{-1}$,
for some essential separating circle $c$ on $S$. 

Suppose, for instance, that $\Psi(D_a) = D_c$. Then $\Psi(D_a)$ is a right
twist. 

Let $b$ be an essential separating circle on $S$. By Proposition
\ref{prop:sepconnected}, there exists a sequence $a_i$, $1 \leq i \leq n$ of
essential separating circles on $S$ such that (i) $a_1 = a$, (ii) $a_i$ is
disjoint from $a_{i + 1}$ for $1 \leq i < n$,  and (iii) $a_n$ is isotopic to
$b$. We may assume, furthermore, that (iv) $a_i$ is not isotopic to $a_{i + 1}$
for $1 \leq i < n$. Let $i$ be an integer such that
$1 \leq i < n$. Then, by the above assumptions, $a_i \cup a_{i + 1}$ is the
boundary  of a compact connected subsurface $R_i$ of $S$ of genus $g_i \geq 1$.
Clearly, by appropriately ``enlarging'' the sequence
$a_i$, $1 \leq i \leq n$, we may assume, furthermore, that (v) $g_i = 1$ for $1
\leq i < n$. 

Let $f_i$ denote the twist about $a_i$, so that $f_1 = D_a$ and $f_n = D_b$. 

Suppose that (i)-(v) hold and let $i$ be an integer such that $1 \leq i < n$.
Then, we may apply Proposition \ref{prop:twoholedtori} to the embedded torus
with two holes $R_i$ in $S$, since both boundary components, $a_i$ and $a_{i +
1}$, of this torus are essential separating circles on $S$. It follows,
therefore, that the restriction of $\Psi$ to $\Gamma (R_i)$ is induced by a
homeomorphism $h : S
\rightarrow S$. Note that $D_i$ and $D_{i + 1}$ are both in $\Gamma (R_i)$.
Hence, $\Psi(D_i) = h_*(D_i)$ and $\Psi(D_{i + 1}) = h_*(D_{i + 1})$. If
$h$ is orientation reversing, then $h_*(D_i)$ and $h_*(D_{i + 1})$ are both left
twists. If $h$ is orientation reversing, then $h_*(D_i)$ and $h_*(D_{i + 1})$
are both right twists. Hence, $\Psi(D_i)$ and $\Psi(D_{i + 1})$ are either both
left twists or both right twists. It follows that the images
$\Psi(D_j)$, $1 \leq j \leq n$ are either all left twists or all right twists. 

Since $\Psi(D_a)$ is a right twist, $\Psi(D_1)$ is a right twist. Hence,
$\Psi(D_n)$ is a right twist. That is to say, $\Psi(D_b)$ is a right twist. In
other words, $\Psi(D_b) = D_d$ for some essential separating circle $d$ on $S$. 

Hence, if $\Psi(D_a) = D_c$ for some essential separating circle $c$ on
$S$, then the result holds with $\epsilon = 1$. Likewise, if
$\Psi(D_a) = D_c^{-1}$ for some essential separating circle $c$ on $S$, then the
result holds with $\epsilon = -1$.
\end{proof} 

Henceforth, we say that an automorphism $\Psi : \mathcal{T} \rightarrow
\mathcal{T}$ is orientation-reversing if $\epsilon = -1$ and
orientation-preserving if $\epsilon = 1$, where $\epsilon$ is the constant in
Proposition \ref{prop:ortype}. In other words, an automorphism
$\Psi : \mathcal{T} \rightarrow \mathcal{T}$ is orientation-reversing
(orientation-preserving) if and only if it sends right twists about essential
separating circles on $S$ to left (right) twists about essential separating
circles on $S$.  We shall refer to the constant $\epsilon$ as the  {\it
orientation type} of $\Psi$. 

We have shown that each automorphism of the Torelli group is either
orientation-preserving or orientation-reversing. In other words, there are no
``hybrid'' automorphisms of the Torelli group, sending some right twists to
right twists and other right twists to left twists. 

\section{Induced automorphisms of the complex of separating curves}
\label{sec:indcomplexsep} 

In this section, we show how any given automorphism of the Torelli group induces
an automorphism of the complex of separating curves on $S$.

By Proposition \ref{prop:ortype}, it follows that each automorphism of
$\mathcal{T}$ induces a permutation of the set of isotopy classes of essential
separating circles on $S$. 

\begin{prop} Let $\Psi : \mathcal{T} \rightarrow \mathcal{T}$ be an
automorphism. Let $\epsilon$ be the orientation type of $\Psi$. Let
$\mathcal{V}_{sep}$ denote the set of isotopy classes of essential separating
circles on $S$. Then: 

\begin{itemize} 
\item There exists a unique function $\Psi_* : \mathcal{V}_{sep}
\rightarrow \mathcal{V}_{sep}$ such that, for each $a$ in
$\mathcal{V}_{sep}$, $\Psi(D_a) = D_b^\epsilon$ where $b =
\Psi_*(a)$
\item $\Psi_* : \mathcal{V}_{sep} \rightarrow \mathcal{V}_{sep}$ is a bijection
\end{itemize} 
\label{prop:indvertsep}
\end{prop} 

Since two essential separating circles have trivial geometric intersection if
and only if the twists about the circles commute, it follows that the induced
map $\Psi_* : \mathcal{V}_{sep} \rightarrow \mathcal{V}_{sep}$ is simplicial. 

\begin{prop} Let $\Psi : \mathcal{T} \rightarrow \mathcal{T}$ be an automorphism
of $\mathcal{T}$. Let
$\Psi_*: \mathcal{V}_{sep} \rightarrow \mathcal{V}_{sep}$ be the unique function
of Proposition \ref{prop:indvertsep}. Then $\Psi_*:
\mathcal{V}_{sep} \rightarrow \mathcal{V}_{sep}$ extends to a simplicial
automorphism $\Psi_* : \mathcal{C}_{sep} \rightarrow \mathcal{C}_{sep}$ of
$\mathcal{C}_{sep}$.
\label{prop:indcomplexsep} \end{prop}

\section{Nonseparating circles}
\label{sec:nonsepcirc} 

In this section, we show how $\Psi$ induces a map of the set
$\mathcal{V}_{nonsep}$ of isotopy classes of nonseparating circles on $S$.  

Suppose that $(a, b)$ is an ordered bounding pair on $S$. That is to say, $a$
and $b$ are elements of $\mathcal{V}_{nonsep}$ represented by disjoint
nonseparating circles, $A$ and $B$, on $S$ such that the complement of $A \cup
B$ in $S$ is disconnected. Note that $A \cup B$ is the boundary of exactly two
embedded surfaces in $S$, $L$ and $R$. Moreover, $L \cap R = A \cup B$. 

By Proposition \ref{prop:presbpmaps}, $\Psi(D_a \circ D_b^{-1}) = D_e \circ
D_f^{-1}$ for some bounding pair $\{e, f\}$ on $S$. Let $\epsilon$ be the
orientation type of $\Psi$. $D_e \circ D_f^{-1} = (D_c \circ D_d^{-1})^\epsilon$
where $\{c, d\} = \{e, f\}$. 

It follows that to each ordered bounding pair $(a,b)$ there exists a unique
ordered bounding pair $(c,d)$ such that $\Psi(D_a \circ D_b^{-1}) = (D_c \circ
D_d^{-1})^\epsilon$. 

\begin{prop} There exists a unique function $\Psi_* : \mathcal{BP} \rightarrow
\mathcal{BP}$ such that if $(a,b)$ is an element of $\mathcal{BP}$ and $(c,d) =
\Psi_*(a,b)$, then $\Psi(D_a \circ D_b^{-1}) = (D_c \circ D_d^{-1})^\epsilon$,
where
$\epsilon$ is the orientation type of $\Psi$.
\label{prop:indbpmaps} \end{prop}

We shall obtain the induced map $\Psi_* : \mathcal{V}_{nonsep} \rightarrow
\mathcal{V}_{nonsep}$ by showing that the first coordinate $c$ of $\Psi_*(a,b)$
does not depend upon the second coordinate $b$ of $(a,b)$.  

To this end, let $a$ be an element of $\mathcal{V}_{nonsep}$ and let 
$\mathcal{BP}_a$ be the set of elements $b$ of $\mathcal{V}_{nonsep}$ such that
$(a,b)$ is an ordered bounding pair on $S$.

Suppose that $b$ and $c$ are distinct elements of $\mathcal{BP}_a$. We say that
$(a,b,c)$ is a {\it $k$-joint based at $a$} if the geometric intersection 
number of $b$ and $c$ is equal to $k$. Since $(a,b)$ and $(a,c)$ are ordered
bounding pairs on $S$, we may orient any representative circles, $A$, $B$, and
$C$, of $a$, $b$, and $c$ so that they are homologous on $S$. Since the oriented
representative circles $B$ and $C$ of $b$ and $c$ are homologous, $i(b,c)$ is
even.

The following result will help us to establish the desired invariance of the
first coordinate of $\Psi_*(a,b)$. 

\begin{prop} Let $(a,b,c)$ be a $k$-joint based at $a$. Let $(e,f) =
\Psi_*(a,b)$ and $(g,h) = \Psi_*(a,c)$. Suppose that $j_1$ and $j_2$ are
distinct elements of $\mathcal{V}_{nonsep}$ such that $i(j_i,a) = i(j_i,b) =
i(j_i,c) = 1$ and $i(j_1,j_2) = 0$. Then $e = g$. 
\label{prop:commonduals} \end{prop}

\begin{proof}
Note that $a$, $b$, $c$, $j_1$ and $j_2$ may be represented by 
transverse nonisotopic nonseparating circles $A$, $B$, $C$, $J_1$, and $J_2$
on $S$ such that $A$ is disjoint from $B$ and $C$; and $J_1$ and $J_2$ intersect
each of $B$ and $C$ in exactly one point. 

Let $P_i$ denote a regular neighborhood of the graph $A \cup J_i \cup B$, $i =
1,2$. Let $Q_i$ denote a regular neighborhood of the graph $A \cup J_i \cup C$,
$i = 1,2$. 

Note that $P_1$ and $Q_1$ are embedded two-holed tori on $S$. Since $(a,b)$ and
$(a,c)$ are bounding pairs, it follows that each of the boundary components of
$P_i$ and $Q_i$ are essential separating circles on $S$. Hence, by Proposition
\ref{prop:twoholedtori}, there exists a homeomorphism $G_i : S \rightarrow S$
such that $\Psi([F]) = [G_i \circ F \circ G_i^{-1}]$ for every homeomorphism $F :
S \rightarrow S$ which is supported on $P_i$ and acts trivially on the first
homology of $S$. Likewise, by Proposition \ref{prop:twoholedtori}, there exists a
homeomorphism $H_i : S \rightarrow S$ such that $\Psi([F]) = [H_i
\circ F \circ H_i^{-1}]$ for every homeomorphism $F : S \rightarrow S$ which is
supported on $Q_i$ and acts trivially on the first homology of $S$. 

Let $R_i$ be a regular neighborhood on $S$ of the graph $A \cup J_i$ such that
$R_i$ is contained in the interiors of both $P_i$ and $Q_i$. Note that $R_i$ is a
one-holed torus. Let $D_i$ be the boundary of $R_i$. Since the genus of $S$ is
not equal to one, $D_i$ is an essential separating circle on $S$. 

Let $U = D_1$ and $u$ be the isotopy class of $D_1$ on $S$. By Propositions
\ref{prop:presseptwists} and \ref{prop:ortype} , $\Psi(D_u) = D_v^\epsilon$ for
some element $v$ of $\mathcal{V}_{sep}$.   

Note that $D_u$ is represented by a twist map $D_U$ supported on $P_1$. Hence,
$\Psi(D_u) = [G_1 \circ D_U \circ G_1^{-1}]$. Note that $[G_1 \circ D_U \circ
G_1^{-1}]$ is equal to $[D_X]^\alpha$ where $X = G_1(U)$ and $\alpha$ is the
orientation type of the homeomorphism $G_1 : S \rightarrow S$. Let $x$ denote
the isotopy class of $X$ on $S$. Then $D_x^\alpha = D_v^\epsilon$. This implies
that $x = v$ and $\alpha = \epsilon$. Hence, $v$ is represented by the circle
$X$ and the homeomorphism $G_1$ has the same orientation type as the automorphism
$\Psi$. 

Likewise, if $Y = H_1(U)$, then $v$ is represented by the circle $Y$ and
the homeomorphism $H_1$ has the same orientation type as the automorphism
$\Psi$. 

Since $v$ is represented by both $X$ and $Y$, $X$ is isotopic to $Y$. Hence, we
may assume that $X = Y$. Let $R_1' = G_1(R_1)$. Note that $R_1'$ is an embedded
one-holed torus with boundary $X$. Since the genus of $S$ is at least $3$, $R_1'$
is the unique embedded one-holed torus in $S$ with boundary $X$. Let $R_1" =
H_1(R_1)$. Note that $R_1"$ is an embedded one-holed torus with boundary $Y$.
Since $X = Y$, $R_1' = R_1"$. 

Note that $D_a \circ D_b^{-1}$ is represented by a bounding pair map $D_A \circ
D_B^{-1}$ where $D_A$ is a twist map supported on $P_1$ and $D_B$ is a twist map
supported on $P_1$. Hence, by the preceding observations, $\Psi(D_a \circ
D_b^{-1}) = (D_p \circ D_q^{-1})^\epsilon$ where $p$ is the isotopy class of
$G_1(A)$ and $y$ is the isotopy class of $G_1(B)$. On the other hand, since
$(e,f) = \Psi_*(a,b)$, $\Psi(D_a \circ D_b^{-1}) = (D_e \circ
D_f^{-1})^\epsilon$. It follows that $e = p$ and $f = q$. Since $e = p$, it
follows that $e$ is represented by the circle $G_1(A)$. Note that $G_1(A)$ is
contained in the interior of the embedded torus $R_1'$. We conclude that $e$ is
represented by the circle $G_1(A)$ contained in the interior of the embedded
torus $R_1'$. 

Likewise, $g$ is represented by the circle $H_1(A)$ and $H_1(A)$ is contained in
the interior of the embedded torus $R_1"$. Since $R_1" = R_1'$, $H_1(A)$ is
contained in the interior of the embedded torus $R_1'$.

It follows that $e$ and $g$ are both represented by circles contained in the
interior of the unique one-holed torus $R_1'$ bounded by $G_1(D_1)$. 

Likewise, $e$ and $g$ are both represented by circles contained in the
interiors of the unique one-hole torus $R_2'$ bounded by $G_2(D_2)$.

Let $d_i$ denote the isotopy class of $D_i$ on $S$. Note that $d_i$ is an
element of $\mathcal{V}_{sep}$. By the preceding considerations, $\Psi_*(d_i) =
e_i$, where $e_i$ is the isotopy class of $G_i(D_i)$, $i = 1, 2$. 

Suppose that $i(e_1, e_2) = 0$. By Proposition \ref{prop:indcomplexsep}, it
follows that $i(d_1,d_2) = 0$. Hence, we may assume, by an isotopy of $D_2$, that
$D_1$ and $D_2$ are disjoint. Since the one-holed tori, $R_1$ and $R_2$, both
contain $A$, it follows that $D_1$ is isotopic to $D_2$. Hence, by a further
isotopy of $D_2$, we may assume that $D_1 = D_2$ and, hence, $R_1 = R_2$. 
Since $J_1$ and $J_2$ are disjoint nonseparating circles, $i(j_1,j_2) = 0$.
Hence, $J_1$ and $J_2$ are isotopic to disjoint nonseparating circles in the
one-holed torus $R_1$. Since any two disjoint nonseparating circles in a
one-holed torus are isotopic, it follows that $j_1 = j_2$. This contradicts our
assumptions. Hence, $i(e_1,e_2)$ is not equal to $0$. 

By an isotopy of $G_2(D_2)$ we may assume that $G_1(D_1)$ and $G_2(D_2)$ are
transverse essential separating circles with minimal intersection. 

By the preceding observations, $i(e,e_i) = i(g,e_i) = 0$, $i = 1,2$. Hence, we
may represent $e$ and $g$ by nonseparating circles $E$ and $G$ on $S$ such that
$E$ and $G$ are each disjoint from both $G_1(D_1)$ and $G_2(D_2)$. 

It follows that $E$ and $G$ are both contained in the interior of the
unique one-holed torus, $R_1'$, with boundary $G_1(D_1)$, and both miss the
intersection of $G_2(D_2)$ with $R_1'$. SInce $i(e_1,e_2)$ is not equal to $0$,
the intersection of $G_2(D_2)$ with $R_1'$ consists of at least one essential
properly embedded arc, $K$. 

Since any two nonseparating circles in a given one-holed torus missing a given
essential properly embedded arc in the given one-holed torus are isotopic in
the given one-holed torus, it follows that $E$ is isotopic to $G$ on $S$. 

Hence, $e= g$. 
\end{proof}

\begin{prop} Let $(a,b,c)$ be a $0$-joint. Let $(e,f) = \Psi_*(a,b)$ and $(g,h) =
\Psi_*(a,c)$. Then $e = g$.
\label{prop:zerojoints} \end{prop}  

\begin{proof} Note that $a$, $b$, and $c$ may be represented by disjoint
nonisotopic nonseparating circles $A$, $B$, and $C$ on $S$. Since $(a,c)$ is a
bounding pair on $S$, $A \cup C$ is the boundary of exactly two embedded surfaces
on $S$, $L$ and $R$. Note that $L$ and $R$ are both connected surfaces and $L
\cap R = A \cup C$. Since $B$ is disjoint from $A$ and $C$, $B$ is contained in
one of the two components of the complement of $A \cup C$ in $S$. That is to
say, either $B$ is contained in the interior of $L$ or $B$ is contained in the
interior of $R$. We may assume that $B$ is contained in the interior of $L$. 

Since $(a,b)$ is a bounding pair on $S$, it follows that $L$ is the union of two
embedded surfaces, $P$ and $Q$, on $S$ such that $A \cup B$ is the boundary of
$P$, $B \cup C$ is the boundary of $Q$, and $P \cap Q = B$. 

Note that the embedded surfaces, $P$, $Q$, and $R$, have disjoint interiors, $P
\cap Q = B$, $Q \cap R = C$, $R \cap P = A$, $P \cap Q \cap R$ is empty, and $P
\cup  Q \cup R = S$. Moreover, $P$, $Q$, and $R$ each have exactly two boundary
components;
$\partial P = A \cup B$, $\partial Q = B \cup C$, and $\partial R = C \cup A$.
Since $a$,
$b$, and $c$ are distinct isotopy classes, $P$, $Q$, and $R$ each have positive
genus.

Since $B$ is not isotopic to $C$, the genus of $Q$ is positive. Hence, there
exists a pair of disjoint, nonseparating circles $D_1$ and $D_2$ on $S$ such
that $Q$ is the union of two embedded connected surfaces, $Q_L$ and $Q_R$, on
$S$, where the interiors of $Q_L$ and $Q_R$ are disjoint; the boundary of $Q_L$
is equal to $B \cup J_1 \cup J_2$; and the boundary of $Q_R$ is equal to $J_1
\cup J_2 \cup C$. 

Choose distinct points, $p_1$ and $p_2$, on $A$; distinct points, $q_1$ and
$q_2$ on $B$; a point $r_1$ on $D_1$; a point $r_2$ on $D_2$; and distinct
points, $s_1$ and $s_2$ on $C$.  Choose disjoint properly embedded arcs, $P_1$
and $P_2$, on $L$ such that $P_1$ joins $p_1$ to $q_1$, and $P_2$ joins $p_2$ to
$q_2$. Choose disjoint properly embedded arcs, $M_1$ and $M_2$, on $Q_L$ such
that $M_1$ joins $q_1$ to $r_1$; and $M_2$ joins $q_2$ to $r_2$. Choose 
disjoint properly embedded arcs, $N_1$ and $N_2$, on $Q_R$ such that $N_1$ joins
$r_1$ to $s_1$; and $N_2$ joins $r_2$ to $s_2$.  Choose disjoint properly
embedded arcs, $R_1$ and $R_2$, on $R$ such that $R_1$ joins $s_1$ to $p_1$ and
$R_2$ joins $s_2$ to $p_2$. 

Let $J_i = P_i \cup M_i \cup N_i \cup R_i$, $i = 1, 2$. Note that $J_1$ and $J_2$
are disjoint circles on $S$; each of these two circles is transverse to $A$, $B$,
$D_1$, $D_2$, and $C$; and each of these two circles intersects each of
the circles, $A$, $B$, and $C$, in exactly one point. It follows that each of
these two circles are nonseparating circles on $S$. 

Let $j_i$ denote the isotopy class of $J_i$ on $S$, $i = 1, 2$. 

Note that each of the circles, $J_i$, is transverse to $D_1$. $J_1$
intersects $D_1$ in exactly one point, whereas $J_2$ is disjoint from $D_1$.
It follows that the geometric intersection of $J_1$ with $D_1$ is equal to $1$,
whereas the geometric intersection of $J_2$ with $D_1$ is equal to $0$. 

It follows that $j_1$ and $j_2$ are distinct elements of
$\mathcal{V}_{nonsep}$ with $i(j_i, a) = i(j_i,b) = i(j_i,c) = 1$, $i = 1,2$, 
and $i(j_1,j_2) = 0$. Hence, by Proposition \ref{prop:commonduals},
$e = g$. 
\end{proof} 

\begin{rem} If the genus of $S$ is at least $5$, it can be shown
that
$\mathcal{BP}_a$ is the vertex set of a connected subcomplex of $\mathcal{C}$.
In other words, any two elements $b$ and $c$ of $\mathcal{BP}_a$ are connected
by a sequence $b_1 = b,...,b_n = c$ of elements of $\mathcal{BP_a}$ such that
$(a,b_i,b_{i+1})$ is a $0-joint$ based at $a$, $1 \leq i < n$. Hence,
Proposition \ref{prop:zerojoints} is sufficient for the purpose of obtaining the
desired induced map $\Psi_* : \mathcal{V}_{nonsep} \rightarrow
\mathcal{V}_{nonsep}$ when $g > 4$. 

One of the main reasons why this result is so much more difficult to achieve
when $g$ is equal to $3$, is that this connectivity result is no longer true
when $g = 3$. (We do not know whether it is true when $g = 4$.) Indeed, when $g
= 3$, there are no edges in $\mathcal{C}$ with both vertices in $\mathcal{BP}_a$.
\label{rem:needtwojoints} \end{rem} 

One of the key ideas for obtaining the main result of this paper for $g > 2$ is
to include with the edges of $\mathcal{C}$ which join elements of
$\mathcal{BP}_a$ more general connections between elements of $\mathcal{BP}_a$
which leave invariant the first coordinate of $\Psi_*(a,b)$ and provide paths
between any two elements of $\mathcal{BP}_a$. 

The first such connection is provided by our next result.

\begin{prop} Let $(a,b,c)$ be a $2$-joint based at $a$. Let $(e,f) =
\Psi_*(a,b)$ and 
$(g,h) = \Psi_*(a,c)$. Then $e = g$.
\label{prop:twojoints} \end{prop} 

\begin{proof} We may represent $a$, $b$, and $c$ by circles, $A$, $B$, and $C$,
on $S$ such that $A$ is disjoint from $B$ and $C$; $B$ is transverse to $C$; and
$B$ and $C$ intersect in exactly two points, $x$ and $y$.

Since $(a,b)$ and $(a,c)$ are bounding pairs, we may orient $A$, $B$, and $C$ so
that the the oriented circles, $A$, $B$, and $C$ are homologous on $S$. Since
$B$ and $C$ are homologous, the algebraic (i.e. homological) intersection of $B$
and $C$ is equal to $0$. It follows that the signs of intersection of $B$ with
$C$ at $x$ and $y$ are opposite.  We assume that the sign of intersection of $B$
with $C$ at $x$ is positive and the sign of intersection of $B$ with $C$ at $y$
is negative.

The above considerations imply that a regular neighborhood $P$ on $S$ of the
graph $B \cup C$ is a four-holed sphere. We may assume that $A$ is contained in
the complement of $P$. 

Since $B$ and $C$ are transverse circles on $S$ with geometric intersection $2$,
it follows that each of the four boundary components of $P$ is an essential
circle on $S$. Indeed, suppose that one of the four boundary components of $P$
bounds a disc $D$ on $P$. Then $P \cup D$ is an embedded pair of pants with $B$
and $C$ in the interior of $P \cup D$. Note that any two circles in a pair of
pants have geometric intersection $0$. Indeed, any circle in a pair of pants on
$S$ is isotopic to one of the three boundary components on the pair of pants
\cite{flp}. Hence, $B$ and $C$ have geometric intersection $0$. This is a
contradiction. Hence, each of the four boundary components of $P$ is an
essential circle on $S$. 

The two points, $x$ and $y$, divide the oriented circle
$B$ into two oriented arcs, $B_1$ and $B_2$, where $B_1$ begins at $x$ and ends
at $y$, $B_2$ begins at $y$ and ends at $x$, $B_1 \cup B_2 = B$, and $B_1 \cap
B_2 = \{x, y\}$. Likewise, the two points, $x$ and $y$ divide the oriented circle
$C$ into two oriented arcs, $C_1$ and $C_2$, where $C_1$ begins at $x$ and ends
at $y$,
$C_2$ begins at $y$ and ends at $x$, $C_1 \cup C_2 = C$, and $C_1 \cap C_2 = \{x,
y\}$. 

Let $i$ and $j$ be integers with $1 \leq i, j \leq 2$. Note that $B_i \cup C_j$
is a circle in the interior of $P$. There is a unique embedded annulus
$A_{(i,j)}$ in $P$ such that the boundary of $A_{(i,j)}$ consists of $B_i \cup
C_j$ and a component
$D_{(i,j)}$ of the boundary of $P$. $P$ is equal to the union of the four annuli,
$A_{(i,j)}$; these annuli have disjoint interiors; and these annuli meet exactly
along those edges of the graph $B \cup C$ which they have in common.

Orient the boundary of $P$ so that $P$ is on the right of its boundary. With this
orientation, the cycle $D_{(1,1)}$ is homologous to $B_1 - C_1$; the cycle
$D_{(1,2)}$ is homologous to $- C_2 - B_1$; the cycle $D_{(2,2)}$ is homologous
to $- B_2 + C_2$; the cycle $D_{(2,1)}$ is homologous to $C_1 + B_2$; the cycle
$B$ is homologous to $B_1 + B_2$; and the cycle $C$ is homologous to $C_1 +
C_2$. By our assumptions, $B$ is homologous to $C$. This implies that $B_1 +
B_2$ is homologous to $C_1 + C_2$; and $B_1 - C_1$ is homologous to $- B_2 +
C_2$. Hence, the cycles $D_{(1,1)}$ and $D_{(2,2)}$ are homologous.  

Note that the cycle $D_{(1,1)} + D_{(1,2)} + D_{(2,2)} + D_{(2,1)}$ is
homologous to $0$, since it represents the boundary of the oriented surface $P$.
This implies that the cycle
$D_{(1,2)} + D_{(2,1)}$ is homologous to $D_{(1,1)} + D_{(2,2)}$. On the other
hand, since $D_{(1,1)}$ is homologous to $D_{(2,2)}$, $D_{(1,1)} + D_{(2,2)}$ is
homologous to
$2D_{(1,1)}$. Hence, $D_{(1,2)} + D_{(2,1)}$ is homologous to $2D_{(1,1)}$. 

Note that the cycle $D_{(1,2)} + D_{(2,1)}$ is homologous to an oriented circle
$E$ contained in the interior of $P$. $E$ is obtained by tubing the oriented
circle $D_{(1,2)}$ to the oriented circle $D_{(2,1)}$ along a properly embedded
arc $J$ in $P$ joining a point on the boundary component $D_{(1,2)}$ of $P$ to a
point on the boundary component
$D_{(2,1)}$ of $P$. 

Suppose that $D_{(1,1)}$ is not homologous to $0$. Then, since the homology of an
oriented surface is torsion free, $2D_{(1,1)}$ is not homologous to $0$. Hence,
$E$ is not homologous to $0$. This implies that $E$ is a nonseparating circle on
$S$. It follows that the homology class of $E$ is not divisible by $2$. Indeed,
$E$ is not divisible by any integer greater than $1$. On the other hand, $E$ is
homologous to $2D_{(1,1)}$. This is a contradiction. Hence, $D_{(1,1)}$ is
homologous to $0$. 

Note that $B$ is homologous to $D_{(1,1)} + D_{(2,1)}$. Since $D_{(1,1)}$ is
homologous to $0$, it follows that $B$ is homologous to $D_{(2,1)}$. Since $A$
is homologous to $B$, it follows that $A$ is homologous to $D_{(2,1)}$. 

Since $A$ is nonseparating, the homology class of $A$ is not $0$. Hence, the
homology class of $D_{(2,1)}$ is not zero. This implies that $D_{(2,1)}$ is an
essential circle on
$S$. 

Let $d$ denote the isotopy class of the essential circle $D_{(2,1)}$ on $S$.
Suppose that
$d$ is not equal to $a$. Since $P$ is disjoint from $A$ and $D_{(2,1)}$ is a
boundary component of $P$, $D_{(2,1)}$ is disjoint from $A$. Since $D_{(2,1)}$
is disjoint from
$A$ and homologous to $A$, it follows that $(a,d)$ is a bounding pair on $S$. 

Note that $D_{(2,1)}$ is also disjoint from $B$ and $C$, since $P$ is a regular
neighborhood of $B \cup C$. Since $i(c,d) = 0$ and $i(c,b) = 2$, $d$ is not
equal to $b$. Hence, $(a,b,d)$ is a $0$-joint based at $a$. Hence, by Proposition
\ref{prop:zerojoints}, the first coordinate of $\Psi_*(a,b)$ is equal to the
first coordinate of $\Psi_*(a,d)$. 

Similarly, $(a,d,c)$ is a $0$-joint based at $a$ and, hence, the first
coordinate of
$\Psi_*(a,d)$ is equal to the first coordinate of $\Psi_*(a,c)$. 

Hence, the first coordinate of $\Psi_*(a,b)$ is equal to the first coordinate of
$\Psi_*(a,c)$. In other words, $e = g$, as desired. 

Hence, we may assume that $d$ is equal to $a$. It follows that $A$ is isotopic to
$D_{(2,1)}$. Since $A$ and $D_{(2,1)}$ are disjoint essential circles on $S$ and
the genus of $S$ is not equal to $1$, it follows that $A \cup D_{(2,1)}$ is the
boundary of exactly one embedded annulus $L$ on $S$. Since $D_{(2,1)}$ is a
boundary component of $P$ and $A$ is in the complement of $P$, it follows that
$L \cap P = D_{(2,1)}$. 

Note that $P \cup L$ is an embedded four-holed sphere on $S$ with boundary
$D_{(1,1)} \cup D_{(1,2)} \cup D_{(2,2)} \cup A$. 

Since $D_{(1,2)} + D_{(2,1)}$ is homologous to $E$ and $E$ is homologous to $0$,
it follows that $D_{(1,2)}$ is homologous to $-D_{(2,1)}$. Since $D_{(2,1)}$ is
homologous to $A$, it follows that $D_{(1,2)}$ is homologous to $-A$. As before,
for $A$ and $D_{(2,1)}$, it follows that $A \cup D_{(1,2)}$ is the boundary of
exactly one embedded annulus $R$ on $S$. 

Since $A$ is a nonseparating circle on $S$, it follows that $A \cup D_{1,2)}$ is
the boundary of exactly two embedded surfaces on $S$. The annulus $R$ is one of
these two surfaces. Let $Q$ be the other of these two surfaces. Note that $R$
and $Q$ are both connected. Hence, either $P \cup L$ is contained in $R$ or $P
\cup L$ is contained in
$Q$. 

Suppose that $P \cup L$ is contained in $R$. Then the essential circles $B$ and
$C$ on
$S$ are both contained in the interior of the annulus $R$. Since $i(b,c) = 2$,
this is a contradiction. Indeed, any two circles in an annulus have zero
geometric intersection number. It follows that $P \cup L$ is contained in $Q$. 

This implies that the four-holed sphere $P \cup L$ and the annulus $R$ have
disjoint interiors and $(P \cup L) \cap R = A \cup D_{(2,1)}$. It follows that
$P \cup L \cup R$ is an embedded two-holed torus, $Z$, on $S$ with boundary
components $D_{(1,1)}$ and
$D_{(2,2)}$. Since $D_{(1,1)}$ and $D_{(2,2)}$ are boundary components of $P$,
they are essential circles on $S$. Hence, we may apply Proposition
\ref{prop:twoholedtori} to the embedded two-holed torus $Z$ on $S$. 

It follows that there is a homeomorphism $H : S \rightarrow S$ such that
$\Psi([F]) = [H
\circ F \circ H^{-1}]$ for every homeomorphism $F : S \rightarrow S$ such that
$F$ is supported on $Z$ and $F$ acts trivially on the first homology of $S$.

Let $U$ be a circle in the interior of $Z$ such that $U$ is isotopic to the
boundary component $D_{(1,1)}$ of $Z$. Let $u$ denote the isotopy class of the
essential separating circle $U$ on $S$. Note that $u$ is an element of
$\mathcal{V}_{sep}$. Let $v = \Psi_*(u)$. Then $\Psi(D_u) = D_v^\epsilon$. 

Note that the element $D_u$ of $\mathcal{T}$ is represented by a twist map $D_U$
supported on $Z$. Hence, $\Psi(D_u) = [H \circ D_U \circ H^{-1}]$. Note that $[H
\circ D_U \circ H^{-1}]$ is equal to $D_w^\alpha$ where $w$ is the isotopy class
of the essential circle $H(U)$ and $\alpha$ is the orientation type of the
homeomorphism $H : S \rightarrow S$. Then $D_w^\alpha = D_v^\epsilon$. This
implies that $\alpha = \epsilon$. Hence, the homeomorphism $H$ has the same
orientation type as the automorphism $\Psi$. 

Note that $A$, $B$, and $C$ are contained in the interior of $Z$. Hence, the
twists, $D_a$, $D_b$, and $D_c$, are represented by twist maps, $D_A$, $D_B$,
and $D_C$, supported on $Z$. It follows that $\Psi(D_a \circ D_b^{-1}) = (D_p
\circ D_q)^\alpha$ where $p$ is the isotopy class of $H(A)$ and $q$ is the
isotopy class of $H(B)$. Furthermore, $\Psi(D_a \circ D_c^{-1}) = (D_p \circ
D_r^{-1})^\alpha$ where $r$ is the isotopy class of $H(C)$. 

Hence, $(D_p \circ D_q^{-1})^\alpha = (D_e \circ D_f^{-1})^\epsilon$ and 
$(D_p \circ D_r^{-1})^\alpha = (D_g \circ D_h^{-1})^\epsilon$. Since $\alpha =
\epsilon$, it follows that $e = p$, $f = q$, $g = p$, and $h = r$. Hence, $e =
g$. \end{proof}

\begin{rem} If the genus of $S$ is at least $4$, it can be shown that if
$(a,b,c)$ is a $2$-joint based at $a$, then $b$ and $c$ are connected by a
sequence $b_1 = b,...,b_n = c$ of elements of $\mathcal{BP_a}$ such that
$(a,b_i,b_{i+1})$ is a $0-joint$ based at $a$. Part of the proof of this fact is
implicit in the proof of Proposition \ref{prop:twojoints}. Following this proof,
it is clear that it remains only to consider the situation when the
representative circles $A$, $B$ and $C$ are contained in a two-holed torus $Z$
as in the last part of this proof.
\label{rem:remtwojoints} \end{rem}

In addition to the connections provided by $0$-joints and $2$-joints based at
$a$, we shall use one more type of connection. This third type of connection is
provided by our next result. 

\begin{prop} Let $(a,b,c)$ be a $4$-joint based at $a$.  Let $(e,f) =
\Psi_*(a,b)$ and $(g,h) = \Psi_*(a,c)$. Then $e = g$.
\label{prop:fourjoints} \end{prop} 

\begin{proof} We may represent the isotopy classes $a$, $b$, and $c$ in
$\mathcal{V}_{nonsep}$ be nonseparating circles, $A$, $B$, and $C$, on $S$ such
that $A$ is disjoint from $B$ and $C$; $B$ and $C$ are transverse; and $B$ and
$C$ intersect in exactly $4$ points. 

Since $(a,b)$ is a bounding pair on $S$, $A \cup B$ is the boundary of exactly
two embedded surfaces $L$ and $R$ on $S$. Note that the interiors of $L$
and $R$ are disjoint and $L \cap R = A \cup B$. Orient the circle $A$ and
$B$ so that the boundary of the oriented surface $L$ is equal to the difference
$A - B$ of the cycles $A$ and $B$. Note that, with the chosen orientations, $L$
is on the left of $A$ and $B$, and $R$ is on the right of $A$ and $B$. Note that
the cycles $A$ and $B$ are homologous on $S$. 

Since $(a,c)$ is a bounding pair on $S$, we may orient the circle $C$ so that $C$
is homologous on $S$ to $A$. Since $A \cup B$ is transverse to $C$, and $A \cup
B$ intersects $C$ in exactly four points, and all four of these points lie on
$B$, it follows that $C \cap R$ is a disjoint union of two properly embedded
arcs, $C_1$ and $C_3$, in $R$; and $C \cap L$ is equal to a disjoint union of two
properly embedded arcs, $C_2$ and
$C_4$, in $L$. We may label these arcs so that the oriented arc $C_1$ joins a
point $w$ on $B$ to a point $x$ on $B$; the oriented arc $C_2$ joins the point
$x$ on $B$ to a point $y$ on $B$; the oriented arc $C_3$ joins the point $y$ on
$B$ to a point $z$ on $B$; and the oriented arc $C_4$ joins the point $z$ on
$B$ to the point $w$ on $B$.

Note that it follows, from the above considerations, that, as we travel along the
oriented circle $C$, the signs of intersection of $B$ with $C$ alternate.
Likewise, by interchanging the roles of $B$ and $C$ in the previous argument, it
follows that the signs of inersection of $B$ with $C$ alternate as we travel
along $B$. Hence, we may choose the notation $C_i$, $1 \leq i \leq 4$, so that
the sign of intersection of $B$ with $C$ at $w$ is positive. 

It follows that the sign of intersection of $B$ with $C$ at $x$ is negative; the
sign of intersection of $B$ with $C$ at $y$ is positive; and the sign of
intersection of $B$ with $C$ at $z$ is negative. Since the signs of intersection
of $B$ with $C$ alternate as we travel along the oriented circle $B$, it follows
that the points of intersection occur in one of the following two orders as we
travel along $B$, beginning at $w$: (i) $(w, x, y, z)$ or (ii) $(w,z,y,x)$.

Suppose that (i) holds. Note that the points of intersection of $B$ with $C$
divide the oriented circle $B$ into oriented arcs, $B_i$, $1 \leq i \leq 4$. We
may choose the notation $B_i$ so that $B_1$ joints $w$ to $x$; $B_2$ joins $x$
to $y$; $B_3$ joints $y$ to $z$; and $B_4$ joins $z$ to $w$.  

It follows that a regular neighborhood on $S$ of the graph $B \cup C$ is a
six-holed sphere $P$. We may assume that $A$ is contained in the complement of
$P$. $P$ is the union of six embedded annuli, $A_i$, $1 \leq i \leq 6$, on $S$.
These annuli have disjoint interiors on $S$. There is a unique boundary
component of $P$, $D_i$, such that $D_i$ is one of the two boundary components
of $A_i$. The other boundary component of $A_i$ is a cycle in the graph $B \cup
C$. We may choose the notation $A_i$, $1 \leq i \leq 6$, so that the boundary
component of $A_1$ contained in the interior of $P$ is the cycle $C_1 - B_1$;
the boundary component of $A_2$ contained in the interior of $P$ is the cycle
$B_2 - C_2$; the boundary component of $A_3$ contained in the interior of $P$ is
the cycle $C_3 - B_3$; the boundary component of $A_4$ contained in the interior
of $P$ is the cycle $B_1 + C_2 + B_3 + C_4$; the boundary component of $A_5$
contained in the interior of $P$ is the cycle $- B_4 - C_3 - B_2 - C_1$; and the
boundary component of $A_6$ contained in the interior of $P$ is the cycle $B_4 -
C_4$. 

Note that the oriented circle $B$ is homologous to the cycle $B_1 + B_2 + B_3 +
B_4$; and the oriented cycle $C$ is homologous to the cycle $C_1 + C_2 + C_3 +
C_4$. By our assumptions, $B$ is homologous to $C$. It follows that $D_1 - D_2
+ D_3 - D_4$ is homologous to $0$. 

On the other hand, $D_1 + D_2 + D_3 + D_4 + D_5 + D_6$ is homologous to $0$.
Hence, $D_5 + D_6$ is homologous to $-2(D_1 + D_3)$. As in the proof of
Proposition \ref{prop:twojoints}, the fact that $D_5 + D_6$ is homologous to a
circle in $P$, obtained by tubing $D_5$ to $D_6$ along a properly embedded arc
$J$ in $P$, implies that $D_1 + D_3$ and $D_5 + D_6$ are both homologous to $0$. 
Again, the homology class of a nonseparating oriented circle on $S$ is primitive
(i.e. not a proper multiple of any nonzero homology class). Since $D_1 + D_2 +
D_3 + D_4 + D_5 + D_6$ is also homologous to $0$, it follows that $D_2 + D_4$ is
homologous to $0$.

Note that $D_2 + D_4$ is homologous to the cycle $B - D_6$. Since $D_2 + D_4$
is homologous to $0$, it follows that $B$ is homologous to $D_6$. Since $B$ is
nonseparating on $S$, $D_6$ is nonseparating on $S$. 

Since $D_5 + D_6$ is homologous to $0$, $D_5$ is homologous to $-D_6$. Since
$D_6$ is a nonseparating circle on $S$ and $D_5$ is a circle disjoint from
$D_6$, it follows that $D_5 \cup D_6$ is the boundary of exactly two embedded
surfaces in $S$. Note that these two surfaces are each connected, have
disjoint interiors, and intersect exactly along $D_5 \cup D_6$. Their
interiors are the two connected components of the complement of $D_5 \cup
D_6$ in $S$. Note that $P$ is contained in one of these two surfaces. Let $Q$
be the other of these two surfaces. Since $Q$ is connected and the boundary of
$Q$ is equal to $D_5 \cup D_6$, there exists a properly embedded arc $M$ on $Q$
joining a point $u$ on $D_5$ to a point $v$ on $D_6$. Note, furthermore that
there exists a properly embedded arc $N$ on $P$ joining $u$ to $v$ such that
$N$ is disjoint from $C$, transverse to $B$, and intersects $B$ in exactly one
point, $p$. Moreover, we may choose $N$ so that $p$ lies in the interior of
$B_4$. 

It follows that $M \cup N$ is a circle on $S$ such that the geometric
intersection of $M \cup N$ with $B$ is equal to $1$, and the geometric
intersection of $M \cup N$ with $C$ is equal to $0$. Note that we may orient the
circle $M \cup N$ so that the algebraic (i.e. homological) intersection of $M
\cup N$ with $B$ is equal to $1$. Since $B$ is homologous to $C$, this implies
that the algebraic intersection of $M \cup N$ with $C$ is equal to $1$. On
the other hand, since $M \cup N$ is disjoint from $C$, the algebraic intersection
of $M \cup N$ with $C$ is equal to $0$. This is a contradiction. 

Hence (ii) holds. Note that the points of intersection of $B$ with $C$
divide the oriented circle $B$ into oriented arcs, $B_i$, $1 \leq i \leq 4$. We
may choose the notation $B_i$ so that $B_1$ joints $w$ to $z$; $B_2$ joins $z$
to $y$; $B_3$ joints $y$ to $x$; and $B_4$ joins $x$ to $w$.  

It follows that a regular neighborhood on $S$ of the graph $B \cup C$ is a
six-holed sphere $P$. We may assume that $A$ is contained in the complement of
$P$. $P$ is the union of six embedded annuli, $A_i$, $1 \leq i \leq 6$, on $S$.
These annuli have disjoint interiors on $S$. There is a unique boundary
component of $P$, $D_i$, such that $D_i$ is one of the two boundary components
of $A_i$. The other boundary component of $A_i$ is a cycle in the graph $B \cup
C$. We may choose the notation $A_i$, $1 \leq i \leq 6$, so that the boundary
component of $A_1$ contained in the interior of $P$ is the cycle $B_1 + C_4$;
the boundary component of $A_2$ contained in the interior of $P$ is the cycle
$- B_4 - C_1$; the boundary component of $A_3$ contained in the interior of $P$
is the cycle $C_2 + B_3$; the boundary component of $A_4$ contained in the
interior of $P$ is the cycle $- C_4 + B_2 - C_2 + B_4$; the boundary component
of $A_5$ contained in the interior of $P$ is the cycle $C_1 - B_3 + C_3 - B_1$;
and the boundary component of $A_6$ contained in the interior of $P$ is the cycle
$- B_2 - C_3$. 

Note that the oriented circle $B$ is homologous to the cycle $B_1 + B_2 + B_3 +
B_4$; and the oriented cycle $C$ is homologous to the cycle $C_1 + C_2 + C_3 +
C_4$. By our assumptions, $B$ is homologous to $C$. It follows that $D_4 - D_5$
is homologous to $0$. Hence, the cycle $D_4$ is homologous to the cycle $D_5$. 

This implies that $D_4 + D_5$ is homologous to $2D_4$. As in the proof
of Proposition \ref{prop:twojoints}, the fact that $D_4 + D_5$ is homologous to a
circle in $P$, obtained by tubing $D_4$ to $D_5$ along a properly embedded arc
$J$ in $P$, implies that $D_4$ and $D_4 + D_5$ are both homologous to $0$. 
Again, the homology class of a nonseparating oriented circle on $S$ is primitive
(i.e. not a proper multiple of any nonzero homology class). Since $D_4$ and $D_4
+ D_5$ are both homologous to $0$, it follows that $D_5$ is homologous to $0$.

Since $B$ and $C$ are transverse circles on $S$ with minimal intersection, it
follows that each of the boundary components, $D_1$, $D_3$, $D_2$, and $D_6$ are
essential circles on $S$ \cite{flp}. 

Since $D_1 + D_2 + D_3 + D_4 + D_5 + D_6$, $D_4$, and $D_5$ are all
homologous to $0$, $D_1 + D_3$ is homologous to $- (D_2 + D_6)$. It follows that
$2(D_1 + D_3)$ is homologous to the cycle $(B_1 + C_4) + (C_2 + B_3) + (B_4 +
C_1) + (B_2 + C_3)$. Hence, $2(D_1 + D_3)$ is homologous to $B + C$. Thus,
$2(D_1 + D_3)$ is homologous to $2B$. This implies that $D_1 + D_3$ is
homologous to
$B$. Hence, $(D_1 + D_3)$ is homologous to $A$. Since $D_1 + D_3$ is
homologous to $-(D_2 + D_6)$, it follows that $(D_2 + D_6)$ is homologous to
$-A$. 

Suppose that $D_1$ and $D_2$ are both separating circles on $S$. Then $D_3$
is homologous to $A$ and $D_6$ is homologous to $-A$. Suppose that the circle
$D_3$ is not isotopic to $A$. Let $d$ be the isotopy class of $D_3$. Then
$(a,b,d)$ and $(a,c,d)$ are both $0$-joints based at $a$. Hence, by Proposition
\ref{prop:zerojoints}, the first coordinates of $\Psi_*(a,b)$ and $\Psi_*(a,c)$
are both equal to the first coordinate of $\Psi_*(a,d)$. Hence, $e = g$. Hence,
we may assume that $D_3$ is isotopic to $A$. Since the genus of $S$ is not
equal to $1$, it follows that $A \cup D_3$ is the boundary of exactly one
embedded annulus $L$ on $S$. Note that $P$ is in the complement of the
interior of $L$. Hence, the boundary of the oriented annulus $L$ is equal to $A
- D_3$. Likewise, we may assume that $D_6$ is isotopic to $-A$ and conclude
that there is exactly one embedded annulus $R$ of $S$ with boundary $A \cup
D_6$ and the boundary of the oriented annulus $R$ is equal to $- (A + D_6)$. 

It follows that $P \cup L \cup R$ is an embedded two-holed torus $Q$ on $S$; $A
\cup B \cup C$ is contained in the interior of $Q$; and each of the two boundary
components of $Q$, $D_2$ and $D_6$, are essential separating circles on $S$.
Hence, by Proposition ref{prop:twoholedtori}, it follows, as in the proof of
Proposition \ref{prop:twojoints}, that $e = g$. 

Hence, we may assume that $D_1$ and $D_2$ are not both separating circles.
Likewise, we may assume that $D_1$ and $D_6$ are not both separating circles
on $S$; $D_3$ and $D_2$ are not both separating circles on $S$; and $D_3$ and
$D_6$ are not both separating circles on $S$. 

Suppose that $D_1$ is a separating circle on $S$. Then $D_3$ is homologous
to $A$. Again, we may assume that $D_3$ is isotopic to $A$. Hence, $D_3 \cup A$
is the boundary of exactly one embedded annulus $L$ in $S$, and the boundary of
the oriented annulus $L$ is equal to $A - D_3$.

Since $D_1$ is a separating circle on $S$, $D_2$ and $D_6$ are both nonseparating
circles on $S$. Since $D_2 + D_6$ is homologous to $-A$, it follows that $D_2
\cup D_6 \cup A$ is the boundary of exactly two embedded surfaces on $S$. The
interiors of these two surfaces are the two connected components of the
complement of $D_1 \cup D_2 \cup A$ in $S$. One of these two surfaces contains
$P$. Let $R$ be the other of these two surfaces. Note that the boundary of the
oriented surface $R$ is equal to $-(A + D_2 + D_6)$.  

Choose distinct points, $p_1$ and $p_2$, on $A$; distinct points, $q_1$ and
$q_2$ on $D_3$; a point $r_1$ on $D_2$; and a point $r_2$ on $D_6$.  Choose
disjoint properly embedded arcs, $P_1$ and $P_2$, on $L$ such that $P_1$ joins
$p_1$ to $q_1$, and $P_2$ joins $p_2$ to $q_2$. Choose disjoint properly
embedded arcs, $Q_1$ and $Q_2$, on $P$ such that $Q_1$ joins $q_1$ to $r_1$;
$Q_2$ joins $q_2$ to $r_2$; $Q_1$ and $Q_2$ are both transverse to $B$ and $C$;
$Q_1$ intersects $B \cup C$ exactly at $x$; and $Q_2$ intersects $B \cup C$
exactly at $y$. Choose disjoint properly embedded arcs, $R_1$ and $R_2$, on
$R$ such that $R_1$ joins $r_1$ to $p_1$ and $R_2$ joins $r_2$ to $p_2$. 

Let $J_i = P_i \cup Q_i \cup R_i$, $i = 1, 2$. Note that $J_1$ and $J_2$ are
disjoint circles on $S$; each of these two circles is transverse to $A$, $B$,
and $C$; and each of these two circles intersects each of the circles, $A$,
$B$, and $C$, in exactly one point. It follows that each of these two circles
are nonseparating circles on $S$. Let $j_i$ denote the isotopy class of $J_i$ on
$S$, $i = 1, 2$. 

Note that each of the circles, $J_i$, is transverse to $D_2$. $J_1$
intersects $D_2$ in exactly one point, whereas $J_2$ is disjoint from $D_2$.
It follows that the geometric intersection of $J_1$ with $D_2$ is equal to $1$,
whereas the geometric intersection of $J_2$ with $D_2$ is equal to $0$. 

It follows that $j_1$ and $j_2$ are distinct common duals to the $4$-joint
$(a,b,c)$ with $i(j_1,j_2) = 0$. Hence, by Proposition \ref{prop:commonduals},
$e = g$. 

Hence, we may assume that $D_1$ is a nonseparating circle on $S$.
Likewise, we may assume that $D_3$, $D_2$, and $D_6$ are nonseparating circles
on $S$.   Since $D_1 + D_3$ is homologous to $A$, and
$D_1$, $D_2$, and $D_3$ are disjoint nonseparating curves, it follows that $D_1
\cup D_2 \cup A$ is the boundary of exactly two embedded surfaces in $S$. These
two surfaces are each connected, and have disjoint interiors. The interiors of
these two surfaces are the two connected components of the complement of $D_1
\cup D_2 \cup A$ in $S$. One of these two surfaces contains $P$. Let $L$ be the
other of these two surfaces. Note that the boundary of the oriented surface $L$
is equal to $A - (D_1 + D_3)$. 

Likewise, $D_2 \cup D_6 \cup A$ is the boundary of exactly two embedded surfaces
in $S$. These two surfaces are each connected, and have disjoint interiors. The
interiors of these two surfaces are the two connected components of the
complement of $D_2 \cup D_6 \cup A$ in $S$. One of these two surfaces contains
$P$. Let $R$ be the other of these two surfaces. Note that the boundary of the
oriented surface $R$ is equal to $- (A  + D_2 + D_6)$. 

As in the previous case, we may construct two disjoint circles, $J_i$, $i = 1,
2$, such that $J_i$ is transverse to $A$, $B$, and $C$; $J_i$ intersects each of
$A$, $B$ and $C$ in exactly one point; $J_i$ is transverse to $D_2$; $J_1$
intersects $D_2$ in exactly one point; and $J_2$ is disjoint from $D_2$. 

Again, it follows that the isotopy classes, $j_i$, $i = 1, 2$ of $J_i$, $i = 1,
2$ are distinct common duals of the $4$-joint $(a,b,c)$ based at $a$ with
$i(j_1, j_2) = 0$. Hence, as in the previous case, it follows, by Proposition 
ref{prop:commonduals}, that $e = g$. 

Hence, in any case, $e = g$. 
\end{proof}

\begin{rem} If the genus of $S$ is at least $5$, it can be shown that if
$(a,b,c)$ is a $4$-joint based at $a$, then $b$ and $c$ are connected by a
sequence $b_1 = b,...,b_n = c$ of elements of $\mathcal{BP_a}$ such that
$(a,b_i,b_{i+1})$ is a $0-joint$ based at $a$. We do not know whether this fact
is true when the genus of $S$ is $4$.
\label{rem:needfourjoints} \end{rem}

Our next result will allow us to show that the connections provided by
Propositions \ref{prop:zerojoints}, \ref{prop:twojoints}, and
\ref{prop:fourjoints} are sufficient to prove the desired invariance of the
first coordinate of $\Psi_*(a,b)$.

\begin{prop} Let $a$ be an element of $\mathcal{V}_{nonsep}$. Let $b$ and
$c$ be elements of $\mathcal{BP}_a$. Then there exists a sequence $b_i$, $1 \leq
i \leq n$ of elements of $\mathcal{BP}_a$ such that (i) $b_1 = b$, (ii) $b_n =
c$, and (iii) $i(b_i,b_{i + 1}) \leq 4$, $1 \leq i < n$. 
\label{prop:shortjoints} \end{prop}

\begin{proof} Represent $a$, $b$, and $c$ by circles $A$, $B$, and $C$ on $S$
such that $A$ is disjoint from $B$ and $C$, $B$ and $C$ are transverse to each
other, and the number of intersection points of $B$ and $C$ is minimal (i.e.
$\#(B \cap C) = i(b,c)$).  Let $k$ be the number of points of intersection of
$B$ and $C$ (i.e. $k = i(b,c)$). The proof is by induction on $k$. 

The assertion of Proposition \ref{prop:shortjoints} clearly holds if $k = 0$. In
this case, let $n = 2$, $b_1 = b$, and $b_2 = c$. Hence, we shall assume that $k
> 0$. 

Orient the circle $A$. Since $(a,b)$ is a ordered bounding pair $A \cup B$ is the
boundary of two embedded surfaces on $S$, $L_B$ and $R_B$, with $L_B \cap R_B =
A \cup B$. Moreover, we may orient $B$ so that the homology classes of the
oriented circles $A$ and $B$ are equal. Likewise, we may orient $C$ so that the
homology classes of the oriented circles $A$ and $C$ are equal. 

Note that the intersection of $C$ with $R_B$ is a disjoint union of properly
embedded arcs in $R_B$, each joining a point of intersection of $B$ and $C$ to
another point of intersection of $B$ with $C$. Each of these properly embedded
arcs on $R_B$ is either separating or nonseparating on $R_B$. 

Let $J$ be one of the components of $C \cap R_B$. Suppose that $J$ separates
$R_B$. Note that exactly one of the two components of the complement of $J$ in
$R_B$ contains the boundary component $A$ of $R_B$. Denote the closure of this
component by $P$. Note that the boundary of $P$ consists of two circles on $S$,
$A$ and a circle $D$. The points $x$ and $y$ on $D$ divide $D$ into two arcs
joining $x$ to $y$, $D_B$ and $D_C$, where $D_B$ lies on $B$ and $D_C$ lies on
$C$. Note that $D_C = J$. Since $B$ and $C$ are transverse essential circles with
minimal intersection, $D$ is an essential circle on $S$.  

Let $B_P$ be a circle in the interior of $P$ such that $B_P$ and $D$ bound an
embedded annulus in $P$. Since $B_P$ is isotopic on $S$ to the essential circle
$D$ on $S$, $B_P$ is also essential on $S$.

Let $Q$ be the closure of the other component of the complement of $J$ in $R_B$.
Note that $P \cap Q = J$. The boundary of $Q$ consists of a single circle $E$.
The points $x$ and $y$ on $E$ divide $E$ into two arcs joining $x$ to $y$,
$E_B$ and $E_C$, where $E_B$ lies on $B$ and $E_C$ lies on $C$. Note that $E_C =
J$, $D_B \cup E_B = B$, and $D_B \cap E_B = \{ x, y \}$. 

Let $B_Q$ be a circle in the interior of $Q$ such that $B_Q$ and $E$ bound an
annulus in $Q$.  As for $D$ and $B_P$, $E$ and $B_Q$ are both essential circles
on $S$. Since $B_P$ and $Q$ bound an annulus in $Q$, we may orient $B_Q$ and $E$
so that the oriented circles $B_P$ and $Q$ are homologous in $Q$. Since
$E$ bounds $Q$, $E$ is nullhomologous on $S$. Hence, $B_Q$ is nullhomologous on
$S$. 

Note that $B \cup B_P \cup B_Q$ is the boundary of a pair of pants $R$ in $S$.
The arc $J$ is properly embedded in $R$, joining the boundary component $B$ of
$R$ to itself, and separating the other two boundary components of $R$, $B_P$
and $B_Q$, from each other. In familiar terminology, $B_P \cup B_Q$ is the result
of tubing $B$ to itself along the arc $J$. 

Note that $C$ is transverse to the submanifold $B_P \cup B_Q$ and the number of
intersection points of $C$ with $B_P \cup B_Q$ is equal to $k - 2$. This implies
that the number of intersection points with $C$ of each of the essential
circles, $B_P$ and $B_Q$, is less than $k$.  

Since the boundary of $P$ consists of $A$ and $D$, we may orient $D$ so that the
oriented circles $A$ and $D$ are homologous. Since $B_P$ is isotopic to $D$, we
may orient $B_P$ so that the oriented circles $B_P$ and $D$ are homologous.
Hence, the oriented circles $B_P$ and $A$ are homologous. On the other hand,
$B_P$ and $A$ are disjoint. 

Suppose that $B_P$ is not isotopic to $A$ on $S$. Let $b'$ be the isotopy class
of $B_P$. By the above considerations, $\{a,b'\}$ is a bounding pair, $i(b,b') =
0$ and $i(b',c) < k$. Hence, the result follows by applying the inductive
hypothesis to the elements $b'$ and $c$ of $\mathcal{BP}_a$. 

Hence, we may assume that $B_P$ is isotopic to $A$ on $S$. Since $B_P$ and
$A$ are disjoint essential circles, it follows that $B_P$ and $A$ bound an
annulus on $S$. Note that the disjoint essential circles, $B_P$ and $A$ bound
exactly two embedded surfaces on $S$. One of these two surfaces is $P$. Since
$B_P$ is contained in the interior of $R_B$, the other of these two surfaces
contains the surface of positive genus, $L_B$. It follows that this other surface
is not an annulus. Hence, $P$ is an annulus.  

Consider the arc $D_B$ of $D$ joining $x$ to $y$ along $B$. Suppose that $C$
intersects $D_B$ at a point $u$ in the interior of $D_B$. Note that there is
exactly one component $J'$ of $C \cap R_B$ with $u$ as one of its two
endpoints.  Note that $J'$ lies in the annulus $P$ and both endpoints of $J'$ lie
in the interior of the arc $D_B$. These two endpoints of $J'$ are therefore
joined by an arc $F$ embedded in the interior of $D_B$. It follows that
$J'$ separates the annulus $P$ into two components, one of which is a disk and
the other of which is an annulus. Since $B$ and $C$ are transverse essential
circles with minimal intersection, the disc component cannot contain the arc $F$
on its boundary. This implies that the disc component contains the arc $J$ on
its boundary. 

It follows that the unique component $P'$ of the complement of $J'$ in $R_B$
which contains $A$ is an annulus contained in $P$. Moreover, the boundary
component $D'$ of $P$ containing $J'$ consists of $J'$ and an arc $D'_B$ of $B$
joining $u$ to a point $v$ in the interior of $D_B$. Moreover, $D'_B$ is
contained in the interior of $D_B$. Indeed, $D'_B$ is the arc $F$. 

It follows that we may assume, by choosing an ``innermost'' arc $J$, that the
interior of $D_B$ contains no points of $C$. 

Let $K$ be the unique component of $C \cap L_B$ such that $y$ is the initial
endpoint of $K$. Since the interior of $D_B$ contains no points of $C$, the
terminal endpoint of $K$ is not in the interior of $D_B$. Suppose that this
endpoint is $x$. Then $C = J \cup K$ and, hence, $B \cap C = \{ x, y \}$. Hence,
since $B$ and $C$ intersect minimally, $i(b, c) = 2$. Clearly, in this
situation, the result holds. 

Hence, we may assume that the terminal endpoint of $K$ is a point $z$ in the
interior of $E_B$. This implies that there is an arc $K'$ parallel to $K$ joining
a point on $B_P$ to a point on $B_Q$ such that $K'$ is contained in the
complement of $C$, $K'$ is transverse to $B$, and $K'$ intersects $B$ in exactly
two points, $u$ and $v$, where $u$ is a point in the interior of $D_B$ near the
endpoint $y$ of $D_B$ and $v$ is a point in the interior of $E_B$ near $z$. 

Let $B'$ be a circle obtained by tubing $B_P$ to $B_Q$ along $K'$. The above
considerations imply that (i) $B'$ is disjoint from $A$, (ii) $B'$ is homologous
to $B$, (iii) $B'$ and $B$ are transverse and have minimal intersection, (iv)
$B'$ and $B$ intersect in exactly $4$ points, (v) $B'$ and $C$ are transverse
and (vi) the number of points of intersection of $B'$ and $C$ is equal to $k -
2$. 

Let $b'$ be the isotopy class of $B'$. By the above considerations, it follows
that $b'$ is not equal to $a$, $(a,b')$ is a bounding pair, $i(b',b) = 4$ and
$i(b',c) \leq k - 2$. The result follows, in this situation, by applying the
inductive hypothesis to $b'$ and $c$. 

It remains to consider the case where there is at least one component $J$ of $C
\cap R_B$ such that the complement of $J$ in $R_B$ is connected. This case is
handled in a manner similar to the previous case. Again, we tube $B$ to itself
along $J$ producing a disjoint union $B_1 \cup B_2$ of disjoint oriented circles
$B_1$ and $B_2$ whose homology classes add up to the homology class of
$B$. 

In this situation, we assume the notations $B_1$ and $B_2$ are chosen so that
$B_1$ corresponds to the oriented arc $N_1$ of the oriented circle $B$ which
joins the initial endpoint $x$ of the oriented arc $J$ to the terminal endpoint
$y$ of $J$ and $B_2$ corresponds to the oriented arc $N_2$ of $B$ joining $y$ to
$x$.

Again, we let $K$ be the unique component of $C \cap L_B$ having the terminal
endpoint $y$ of $J$ as its initial endpoint. Again, if the terminal endpoint of
$K$ is equal to $x$, then $C = J \cup K$, $i(b,c) = 2$, and the result clearly
holds. 

Hence, we may assume that the terminal endpoint of $K$ is a point $z$ in the
interior of $N_i$ for some integer $i$ in $\{ 1, 2\}$. Suppose that $z$ is in
the interior of $N_1$. Then we may tube $B_2$ to $B_1$ along an arc $K'$ parallel
to $K$, crossing $N_2$ at a point $u$ in the interior of $N_2$ near $y$ and $N_1$
at a point $v$ in the interior of $N_1$ near $z$. This results in a curve $B'$
with which the proof proceeds by induction as before. 

Suppose, on the other hand, that $z$ is in the interior of $N_2$. In this case,
we may tube $B_1$ to $B_2$ along an arc $K'$ parallel to $K$, crossing $N_1$ at a
point $u$ in the interior of $N_1$ near $y$ and $N_2$ at a point $v$ in the
interior of $N_2$ near $z$. This results in a curve $B'$ with which the proof
proceeds by induction as before.
\end{proof}

We are now ready to establish the desired invariance of the first coordinate of
$\Psi_*(a,b)$. 

\begin{prop} Let $a$ be an element of $\mathcal{V}_{nonsep}$.  Suppose that $b$
and
$c$ are distinct elements of $\mathcal{BP}_a$. Let $(e,f) = \Psi_*(a,b)$ and
$(g,h) =
\Psi_*(a,c)$. Then $e = g$. 
\label{prop:firstcoord} \end{prop}

\begin{proof} By Proposition \ref{prop:shortjoints}, there exists a sequence
$b_i$, $1
\leq i \leq n$ of elements of $\mathcal{BP}_a$ such that $b_1 = b$, $b_n = c$,
and 
$i(b_i,b_{i + 1}) \leq 4$, $1 \leq i \leq n$. 

Suppose that $i$ is an integer with $1 \leq i < n$. Let $(e_i,f_i) =
\Psi_*(a,b_i)$. Note that $(e_1,f_1) = (e,f)$ and $(e_n,f_n) = (g,h)$. 

By the above observations, $(a,b_i,b_{i+1})$ is a $k$-joint based at $a$ for
some even integer $k$ with $0 \leq k \leq 4$. Hence, $k$ is an element of
$\{0,2,4\}$. It follows, from Propositions \ref{prop:zerojoints},
\ref{prop:twojoints}, and \ref{prop:fourjoints}, that $e_i = e_{i + 1}$. 

Thus, by induction, $e_1 = e_n$. That is to say, $e = g$.
\end{proof}

The preceding results now provide us with the desired induced map $\Psi_* :
\mathcal{V}_{nonsep} \rightarrow \mathcal{V}_{nonsep}$.

\begin{prop} Let $\Psi : \mathcal{T} \rightarrow \mathcal{T}$ be an automorphism
of $\mathcal{T}$. Let $\epsilon$ be the orientation type of $\Psi$. Let
$\mathcal{V}_{nonsep}$ denote the set of isotopy classes of nonseparating
circles on $S$. Then: 

\begin{itemize} 
\item There exists a unique function $\Psi_* : \mathcal{V}_{nonsep}
\rightarrow \mathcal{V}_{nonsep}$ such that, for each bounding pair $(a,b)$ on
$S$, $\Psi(D_a \circ D_b^{-1}) = (D_c \circ D_d^{-1})^\epsilon$, where
$c = \Psi_*(a)$ and $d = \Psi_*(b)$. 
\item $\Psi_* : \mathcal{V}_{nonsep} \rightarrow \mathcal{V}_{nonsep}$ is a
bijection.
\end{itemize} 
\label{prop:indvertnonsep}
\end{prop} 

\begin{proof} By Proposition \ref{prop:firstcoord}, there exists a well-defined
function
$\Psi_* : \mathcal{V}_{nonsep} \rightarrow \mathcal{V}_{nonsep}$ defined by the
rule
$\Psi_*(a) = c$ if and only if for each element $b$ of $\mathcal{BP}_a$,
$\Psi_*(a,b) = (c,d)$ for some element $d$ in $\mathcal{BP}_d$. 

Let $(a,b)$ be a bounding pair on $S$. Let $(c,d) = \Psi_*(a,b)$. By the
definition of
$\Psi_* : \mathcal{V}_{nonsep} \rightarrow \mathcal{V}_{nonsep}$, $c =
\Psi_*(a)$. On the other hand, by the definition of $\Psi_* : \mathcal{BP}
\rightarrow \mathcal{BP}$ given in Proposition \ref{prop:indbpmaps}, $\Psi(D_a
\circ D_b^{-1}) = (D_c \circ D_d^{-1})^\epsilon$. Since $\Psi : \mathcal{T}
\rightarrow \mathcal{T}$ is a homomorphism, this implies that $\Psi(D_b \circ
D_a^{-1}) = (D_d \circ D_c^{-1})^\epsilon$. Hence, by the definition of $\Psi_*
: \mathcal{BP} \rightarrow
\mathcal{BP}$, $\Psi_*(b,a) = (d,c)$. Hence, by the definition of $\Psi_* :
\mathcal{V}_{nonsep} \rightarrow \mathcal{V}_{nonsep}$, $d = \Psi_*(b)$. 

This proves the existence of a function $\Psi_*$ as described in the first
clause of Proposition \ref{prop:indvertnonsep}. 

Suppose that $\Psi_{\#} : \mathcal{V}_{nonsep} \rightarrow \mathcal{V}_{nonsep}$
is a function such that, for each bounding pair $(a,b)$ on $S$, $\Psi(D_a \circ
D_b^{-1}) = (D_c \circ D_d^{-1})^\epsilon$, where $c = \Psi_{\#}(a)$ and $d =
\Psi_{\#}(b)$. Let $a$ be an element of $\mathcal{V}_{nonsep}$. Let $b$ be an
element of $\mathcal{BP}_a$. Then
$(a,b)$ is a bounding pair on $S$. Hence, $\Psi(D_a \circ D_b^{-1}) = (D_c \circ
D_d^{-1})^\epsilon$ and $\Psi(D_a \circ D_b^{-1}) = (D_e \circ
D_f^{-1})^\epsilon$, where 
$c = \Psi_*(a)$, $d = \Psi_*(b)$, $e = \Psi_{\#}(a)$, and $f = \Psi_{\#}(b)$.
Note that 
$(D_c \circ D_d^{-1})^\epsilon = (D_e \circ D_f^{\#})^\epsilon$. This implies
that $(c,d) = (e,f)$. Hence, $c = e$. That is to say, $\Psi_*(a) =
\Psi_{\#}(a)$. Hence, the two functions $\Psi_* : \mathcal{V}_{nonsep}
\rightarrow \mathcal{V}_{nonsep}$ and 
$\Psi_{\#} : \mathcal{V}_{nonsep} \rightarrow \mathcal{V}_{nonsep}$ are equal. 

This proves the uniqueness of a function $\Psi_*$ as described in the first
clause of Proposition \ref{prop:indvertnonsep}. 

Let $\Theta : \mathcal{T} \rightarrow \mathcal{T}$ be the inverse of the
automorphism
$\Psi: \mathcal{T} \rightarrow \mathcal{T}$. Clearly, by the definition of the
orientation type of an automorphism in Proposition \ref{prop:ortype} and the
definition of the induced map on $\mathcal{V}_{nonsep}$ given in the first
clause of Proposition \ref{prop:indvertnonsep}, the induced maps $\Psi_* :
\mathcal{V}_{nonsep} \rightarrow \mathcal{V}_{nonsep}$ and $\Theta_* :
\mathcal{V}_{nonsep} \rightarrow \mathcal{V}_{nonsep}$ are inverse maps. This
proves the second clause of Proposition \ref{prop:indvertnonsep}.
\end{proof}

\section{Induced automorphism of the complex of curves}
\label{sec:indautcomplex}

In this section, we extend the automorphism $\Psi_* : \mathcal{C}_{sep}
\rightarrow \mathcal{C}_{sep}$ of Proposition \ref{prop:indcomplexsep} to the
entire complex of circles $\mathcal{C}$ of $S$.

Consider the unique function $\Psi_* : \mathcal{V} \rightarrow \mathcal{V}$ whose
restrictions to $\mathcal{V}_{sep}$ and $\mathcal{V}_{nonsep}$ are the functions
$\Psi_* : \mathcal{V}_{sep} \rightarrow \mathcal{V}_{sep}$ and 
$\Psi_* : \mathcal{V}_{nonsep} \rightarrow \mathcal{V}_{nonsep}$ of Propositions
\ref{prop:indvertsep} and \ref{prop:indvertnonsep}. 

We shall show that $\Psi_* : \mathcal{V} \rightarrow \mathcal{V}$ extends to a
simplicial map. This will use the following consequence of Proposition
\ref{prop:twoholedtori}. 

\begin{prop} Let $a$ be an element of $\mathcal{V}_{nonsep}$ and $b$ be an
element of $\mathcal{sep}$. Let $c = \Psi_*(a)$ and $d =
\Psi_*(b)$. Suppose that $a$ is represented by a circle $A$ on $S$ and $b$ is
represented by a circle $B$ on $S$ such that $B$ is the boundary of an embedded
torus $P$ on $S$ and $A$ is contained in the interior of $P$. Then
$c$ is represented by a circle $C$ on $S$ and $d$ is represented by a circle $D$
on $S$ such that $D$ is the boundary of an embedded torus $Q$ on
$S$ and $C$ is contained in the interior of $Q$.
\label{prop:preslollipops} \end{prop} 

We shall now show that $\Psi_* : \mathcal{V} \rightarrow \mathcal{V}$ extends to
a simplicial map $\Psi_* : \mathcal{C} \rightarrow \mathcal{C}$.

\begin{thm} Let $\Psi : \mathcal{T} \rightarrow \mathcal{T}$ be an automorphism
of
$\mathcal{T}$. Let $\Psi_* : \mathcal{V} \rightarrow \mathcal{V}$ be the unique
function whose restrictions to $\mathcal{V}_{sep}$ and $\mathcal{V}_{nonsep}$
are the functions
$\Psi_* : \mathcal{V}_{sep} \rightarrow \mathcal{V}_{sep}$ and 
$\Psi_* : \mathcal{V}_{nonsep} \rightarrow \mathcal{V}_{nonsep}$ of Propositions
\ref{prop:indvertsep} and \ref{prop:indvertnonsep}.  Then 
$\Psi_*: \mathcal{V} \rightarrow \mathcal{V}$ extends to a simplicial
automorphism 
$\Psi_* : \mathcal{C} \rightarrow \mathcal{C}$ of $\mathcal{C}$.
\label{thm:indautcomplex} \end{thm}

\begin{proof} Let $a$ and $b$ be distinct elements of $\mathcal{V}$ such that
$i(a,b) = 0$. Let $c = \Psi_*(a)$ and $d = \Psi_*(b)$. It suffices to show that
$i(c,d) = 0$. This will follow by considering several cases. 

Suppose that $a$ and $b$ are both elements of $\mathcal{V}_{sep}$. Then, by
Proposition \ref{prop:indcomplexsep}, $i(c,d) = 0$. 

Suppose, on the other hand, that $a$ and $b$ are both elements of
$\mathcal{V}_{nonsep}$. If $(a,b)$ is a bounding pair, then by Proposition
\ref{prop:indvertnonsep}, $(c,d)$ is a bounding pair. In particular,
$i(c,d) = 0$. We may assume, therefore, that $(a,b)$ is not a bounding pair.
Since $i(a,b) = 0$, this assumption implies that $a$ is represented by a circle
$A$ on $S$ and $b$ is represented by a circle $B$ on $S$ such that $A$ is
disjoint from $B$ and the complement of $A \cup B$ in $S$ is connected. Hence,
there exists circles $E$ and $F$ on $S$ such that (i)
$A$, $B$, $E$, and $F$ are all disjoint, (ii) $E$ and $F$ bound embedded tori
$P_E$ and $P_F$ on $S$, and (iii) $A$ is contained in the interior of
$P_E$ and $B$ is contained in the interior of $P_F$. 

Since the genus of $S$ is at least $3$, $E$ is not isotopic to $F$. Let $e$ and
$f$ denote the isotopy classes of $E$ and $F$. Since $E$ and
$F$ are separating circles, $e$ and $f$ are elements of
$\mathcal{V}_{sep}$. Since $E$ and $F$ are disjoint, $i(e,f) = 0$. Since
$E$ and $F$ are not isotopic, $e$ is not equal to $f$.   

Let $g = \Psi_*(e)$ and $h = \Psi_*(f)$. By Proposition
\ref{prop:indvertsep}, since $e$ and $f$ are distinct elements of
$\mathcal{V}_{sep}$, 
$g$ and $h$ are distinct elements of $\mathcal{V}_{sep}$. By Proposition
\ref{prop:indcomplexsep}, since $i(e,f) = 0$, we have $i(g,h) = 0$.

By Proposition \ref{prop:preslollipops}, $c$ and $g$ are represented by circles
$C$ and $G$ on $S$ such that $G$ bounds an embedded torus $Q_G$ on
$S$ and $C$ is contained in the interior of $Q_G$. Likewise, $d$ and $h$ are
represented by circles $D$ and $H$ on $S$ such that $H$ bounds an embedded torus
$Q_H$ on $S$ and $D$ is contained in the interior of $Q_H$. 

Since $i(g,h) = 0$, we may assume that $G$ and $H$ are disjoint. Since $g$ is
not equal to $h$, $G$ is not isotopic to $H$. Hence, the interiors of the
embedded tori, $Q_G$ and $Q_H$, bounded by $G$ and $H$ are disjoint. This
implies that $C$ and $D$ are disjoint. Hence, 
$i(c,d) = 0$.   

Thus, if $a$ and $b$ are both elements of $\mathcal{V}_{nonsep}$, then
$i(c,d) = 0$. 

Suppose, finally, that one of $a$ and $b$ is an element of
$\mathcal{V}_{nonsep}$ and one is an element of $\mathcal{V}_{sep}$. We may
assume that $a$ is an element of $\mathcal{V}_{nonsep}$ and $b$ is an element of
$\mathcal{V}_{sep}$. Since $i(a,b) = 0$, $a$ is represented by a circle $A$ on
$S$ and $b$ is represented by a circle $B$ on $S$ such that
$A$ and $B$ are disjoint. Since $b$ is in
$\mathcal{V}_{sep}$, $B$ bounds an embedded surface $P_B$ in $S$ such that
$A$ is contained in the interior of $P_B$. $P_B$ is the closure of the unique
component of the complement of $B$ in $S$ which contains $A$. Let
$k$ be the genus of $P_B$. Since $B$ is an essential circle on $S$,
$k > 0$. If $k = 1$, then by Proposition \ref{prop:preslollipops}, $i(c,d) = 0$.
Hence, we may assume that $k > 1$. It follows that there exists a circle $E$ on
$S$ such that (i) $A$, $B$, and $E$ are disjoint, (ii) $E$ bounds an embedded
torus $P_E$ on $S$ such that $P_E$ is contained in the interior of $P_B$ and $A$
is contained in the interior of $P_E$. Note that
$E$ is not isotopic to $B$. 

Let $e$ denote the isotopy class of $E$ on $S$. Let $f = \Psi_*(e)$. Since
$B$ and $E$ are nonisotopic disjoint essential separating circles on
$S$, it follows, by Propositions \ref{prop:indvertsep} and
\ref{prop:indcomplexsep}, that $d$ and $f$ are distinct elements of
$\mathcal{V}_{sep}$ such that $i(d,f) = 0$. Hence, $d$ and $f$ are represented
by nonisotopic disjoint essential separating circles, $D$ and
$F$, on $S$.  Hence, by Proposition \ref{prop:preslollipops}, it follows that
$c$ is represented by a circle $C$ on $S$ such that $C$ is contained in the
interior of an embedded torus $P_F$ on $S$ with boundary $F$. Since
$D$ and $F$ are disjoint nonisotopic essential separating circles on $S$, it
follows that $D$ is contained in the complement of $P_F$. Since $C$ is contained
in the interior of $P_F$, it follows that $C$ and $D$ are disjoint. Hence,
$i(c,d) = 0$. 
\end{proof}


\begin{thebibliography}{XXXX}

\bibitem[B]{b} J. S. Birman, {\it Braids, Links and Mapping Class Groups},
Annals of Math. Studies, No. 82, Princeton University Press, Princeton, New
Jersey, 1974.

\bibitem[BLM]{blm} J. S. Birman, A. Lubotzky, and J. McCarthy, Abelian and
solvable subgroups of the mapping class group, Duke Math. J. V. 50 (1983),
1107-1120.

\bibitem[DG]{dg} J. L. Dyer and E. K. Grossman, The automorphism groups of the
braid groups, Amer. J. Math. 103 (1981), 1151-1169.

\bibitem[EE]{ee} C. J. Earle and J. Eells, A fibre bundle description of
Teichm\"{u}ller theory, J. Diff.  Geometry, V. 3, No. 1 (1969), 19-43.

\bibitem[F]{f} B. Farb, Automorphisms of the Torelli group, preliminary report - AMS Sectional Meeting,
Ann Arbor,  Michigan, March 1, 2003

\bibitem[FIv]{fiv} B. Farb, and N. V. Ivanov, The Torelli geometry and its applications, 
preprint, arXiv: math.GT/0311123  

\bibitem[FLP]{flp} A. Fathi, F. Laudenbach. and V. Po\'{e}naru, {\it Travaux de
Thurston  sur les surfaces}, S\'eminaire Orsay, Ast\'{e}risque, Vol. 66-67, Soc.
Math. de France, 1979.

\bibitem[I]{i} E. Irmak, Superinjective Simplicial Maps of Complexes of Curves and
Injective Homomorphisms of Subgroups of Mapping Class Groups, preprint, arXiv: math.GT/0211139, 
to appear, Topology  

\bibitem[IIvM]{iivm} E. Irmak, N. V. Ivanov, and J. D. McCarthy, Automorphisms of surface braid
groups,  preprint, arXiv: math.GT/0306069 

\bibitem[Iv1]{iv1} N. V. Ivanov, Complexes of curves and the Teichm\"{u}ller
modular group, Uspekhi Mat. Nauk V. 42, No. 3 (1987), 49-91; English transl.:
Russian Math. Surveys V. 42, No. 3 (1987),  55-107.

\bibitem[Iv2]{iv2} N. V. Ivanov, Automorphisms of Teichm\"{u}ller modular groups,
Lecture Notes in Math., No. 1346, Springer-Verlag, Berlin and New York, 1988,
199-270.

\bibitem[Iv3]{iv3} N. V. Ivanov, {\it Subgroups of Teichm\"{u}ller Modular
Groups},  Translations of Mathematical Monographs, Vol. 115, American Math.
Soc., Providence, Rhode Island,  1992.

\bibitem[Iv4]{iv4} N. V. Ivanov, Permutation representations of braid groups of
surfaces (Russian)  Mat. Sb. 181 (1990), no. 11, 1464--1474; translation in
Math. USSR-Sb. 71 (1992), no. 2, 309--318

\bibitem[IvM]{ivm} N. V. Ivanov and J. D. McCarthy, On injective homomorphisms
between Teichm\"{u}ller modular groups, I,  Invent. math. 135, 425-486 (1999)

\bibitem[K]{k} M. Korkmaz, Automorphisms of complexes of curves on punctured spheres and on
punctured tori. {\it Topology Appl.,} {\bf 95} (1999), 85-111.

\bibitem[M]{m} J. D. McCarthy, Automorphisms of surface mapping class groups. A
recent theorem of N. Ivanov, Invent. Math. V. 84, F. 1 (1986), 49-71.

\bibitem[MV]{mv} J. D. McCarthy and W. R. Vautaw, Automorphisms of the Torelli
Group in genus $>$ 2, preliminary report - AMS Sectional Meeting,
Baton Rouge, Louisiana, March 14, 2003

\bibitem[Me]{me} G. Mess, The Torelli group for genus 2 and 3 surfaces, Topology 31
(1992), 775-790

\bibitem[V1]{v1} W. R. Vautaw, Abelian Subgroups and Automorphisms of the
Torelli Group, Ph. D. Thesis, Michigan State University, 2002 

\bibitem[V2]{v2} W. R. Vautaw, Abelian Subgroups and Automorphisms of the
Torelli Group:Addendum - A Centralizer Lemma, preprint, December 16, 2002

\bibitem[ZVoC]{zvoc} H. Zieschang, E. Vogt, H.-D. Coldeway, {\it Surfaces and
planar discontinuous groups}, Lecture Notes in Math., No. 835, Springer-Verlag,
Berlin, Heidelberg, New York, 1980. 

\end{thebibliography}
\end{document}